\documentclass{article}
\usepackage{amsmath,amsfonts,amssymb,amsthm}
\usepackage{booktabs}
\usepackage{url}
\usepackage{graphicx}
\usepackage[round]{natbib}

\renewcommand{\le}{\leqslant}
\renewcommand{\ge}{\geqslant}

\newcommand{\e}{\mathbb{E}}
\newcommand{\var}{\mathrm{var}}

\newcommand{\dnorm}{\mathcal N}

\newcommand{\sign}{\mathrm{sign}}

\newcommand{\tod}{\stackrel{\mathrm{d}\ }{\to}}

\newcommand{\giv}{\!\mid\!}

\date{November 2019}
\title{Refiltering hypothesis tests to control sign error}
\author{Art B. Owen\\Stanford University}
\begin{document}
\maketitle

\begin{abstract}
A common, though not recommended statistical practice
is to report confidence intervals if and only if they exclude a null
value of $0$.
The resulting filtered confidence intervals generally do not
have their nominal confidence level. More worryingly, in
low  power settings their center points will be much farther
from zero than the true parameter is and
they will frequently lie on the wrong side of zero.
Many confidence intervals are constructed
using an asymptotically Gaussian parameter estimate accompanied
by a weakly consistent estimate of its variance.  In these cases,
we can subject the given confidence interval(s) to a second
filtering step such that the probability of a sign error
is controled.  This refiltering step retains only
those confidence intervals that are sufficiently well
separated from the origin.  It requires no assumptions on
the dependencies among the test statistics.
\end{abstract}

\smallskip
\par\noindent

\section{Introduction}

A serious problem arises for hypothesis tests of low statistical power,
as pointed out by \cite{gelm:2014:power06}.
He considers an extreme case of testing a null
hypothesis that $\theta=0$ at significance level $\alpha=0.05$
with power only $0.06$.  The type II error is then $\beta=0.94$.
In his formulation, reproduced below,
any rejected hypotheses come with a very exaggerated estimate
of the magnitude of $\theta$.  They frequently have the wrong
sign too.

These points were made earlier by~\cite{gelm:tuer:2000}.
Getting a statistically significant result with
the incorrect sign is sometimes called a directional
error.  It is also sometimes called a
type III error, though that term is also used
to include tests that answer what is in some sense the `wrong' question.
Sign errors have potentially more serious consequences than the type I error
that statisticians focus on.
Sign errors have been discussed in the statistical and related
literatures for a long time. See for instance,
\cite{shaffer2002multiplicity} and \cite{jones2000sensible} and
references therein going back to the 1950s and 1960s.

In this article, we mitigate the sign error problem by
inspecting the confidence interval for $\theta$.
We suppose that somebody has run
one or more hypothesis tests and reports their confidence
intervals, but only those confidence intervals that do
not include the origin.  We will see that when a confidence
interval is well separated from the origin, then we
can be confident that the sign of its center point is correct.
In a refiltering step, we can declare that $\sign(\theta)=\sign(\hat\theta)$
only in those well separated intervals.

Given a confidence interval,
we do not know what power the corresponding test had.
As mentioned above, we can reason that for a confidence
interval that is well separated from the origin we have
more cause to believe that the sign is correct, than
we do for a confidence interval that just barely excludes the origin.
When the tests and confidence intervals are based on
an asymptotically normal test statistic $\hat\theta$
with a suitable  estimate of $\var(\hat\theta)$ then
we can say how far the confidence interval must
be from the origin, in units of its own length, for
us to have some confidence that $\sign(\theta)=\sign(\hat\theta)$.
That in turn leads us to a way to get a $p$-value for
the sign of $\theta$ that takes into account
the fact that we were only presented with confidence intervals
that excluded~$0$.

The remainder of this paper is as follows.
Section~\ref{sec:model} introduces our model formulation.
We consider power for a test of a scalar $\theta$ for which
an asymptotically normal estimate $\hat\theta$ is available
along with a consistent estimator of $\var(\hat\theta)$.
When the confidence interval for $\theta$ is well separated
from the origin, as described in Section~\ref{sec:ciwidth},
then we can be more confident that $\sign(\theta)=\sign(\hat\theta)$.
Equivalently, that confidence can come from repeating
the test at a more stringent level. The result is a criterion
for deciding when to declare the sign of $\theta$ subject
to an upper bound on the probability of making an erroneous
declaration.
Section~\ref{sec:power} computes and displays some operating
characteristics of this testing procedure.
Section~\ref{sec:signp} presents some `sign $p$-values'.
These are the smallest upper bounds on the probability of a sign
error given the selected confidence interval.
Section~\ref{sec:onetail} briefly considers one-sided tests.
Section~\ref{sec:conc} presents some conclusions.
There is an appendix with some R code
to handle some of the numerical tasks in this paper.

To conclude this section, we note that
there are more sophisticated ways to handle multiple confidence
intervals than simply reporting the ones that do not
contain the origin.
For instance, \cite{benj:yeku:2005} propose a confidence interval
version of the well-known procedure of \cite{benj:hoch:1995}.
Such methods are not standard and so it remains
worthwhile to have a means of adjusting intervals
that were selected by non-coverage of zero.
Also, methods designed to handle a multiplicity of
confidence intervals typically assume independence
or some other strong condition on the corresponding
test statistics.  Independence is difficult to achieve
in high throughput settings.  For instance in a genomic
setting measuring a Gaussian response for
$n$ subjects on $p$ genes, there can be at most
$p-1$ independent linear test statistics, while those
settings commonly have $p\gg n$.
The refiltering proposal presented
here does not require independence.
For notions of the false sign proportion in a set of
decisions and the corresponding
false sign rate, see \cite{wein:ramd:2019} and articles
cited there. \cite{step:2016} presents a local false
sign rate to parallel the local false discovery
rate of \cite{efro:2005}.

\section{Model formulation}\label{sec:model}

In our setting there is a real-valued parameter $\theta$
and the null hypothesis is  $H_0:\theta=0$.
We wish to test this hypothesis.
In many contexts, $\theta=0$ means that some
phenomenon is unimportant. For instance, $\theta$
might be the coefficient of $x$  in a regression model
for $y$ and $\theta=0$ then means that $x$ does not
affect the expected value of $y$.

While $|\theta|$ could well be small,
it is a priori almost certain that $\theta$ is not exactly
zero. Even if there is a scientific reason for $\theta=0$
to be possible, we are extremely unlikely to have perfect instruments,
and so the $\theta$ governing our data will not be zero.
Given that $\theta\ne0$, it becomes interesting
to justify doing the test.  A potential explanation is that
if we reject $H_0$, then we have learned the sign of $\theta$
about which there may have been a reasonable doubt.
\cite{jones2000sensible} formulate the problem as one
of testing $\theta\ge0$ versus $\theta\le0$ and show
that the sign error probability is below $\alpha/2$.
Unfortunately, things are more difficult if the interval
we get has been selected for non-containment of the origin.
As we will see, in the low power setting, rejecting
$H_0$ can leave substantial  doubt about $\sign(\theta)$.
Another explanation is that while we may be confident
that $\theta\ne 0$, we may still have to furnish evidence
of that to other people.  For the present purposes, we work
with $\theta\ne0$ to measure the extent to which
the significance filter gives large but finite values of $|\hat\theta/\theta|$.

To formulate a model for this selective process,
we assume that there is an
estimate $\hat\theta\sim\dnorm(\theta,\tau_0^2)$
for some $\tau_0>0$.
The user has an estimator $s_0^2$ of $\tau_0^2$
with $s_0^2/\tau_0^2\tod 1$.
For instance, it might hold that $s_0^2\sim\tau_0^2\chi^2_{(\nu)}/\nu$
with $\nu\to\infty$.
Then $\Pr( |s_0/\tau_0-1|>\epsilon)\to0$ for all $\epsilon>0$.
To focus on essentials, we will work as if $s_0=\tau_0$.
The test of $H_0$ is
based on $Z=\hat\theta/s_0$.  It is a two-tailed test
at level $\alpha$ that rejects $H_0$ when $\hat\theta^2/s_0^2\ge \chi^{2,1-\alpha}_{(1)}$.
It has the same asymptotic rejection probability
as the test that rejects when $\hat\theta^2/\tau_0^2\ge\chi^{2,1-\alpha}_{(1)}$.

Because $\theta\ne 0$, we may change units
in our analysis to make $|\theta|=1$.
Then
$$
\frac{\hat\theta}{\tau_0} =
\frac{\hat\theta/|\theta|}{\tau} \sim\dnorm\Bigl(\frac{\sign(\theta)}{\tau},1\Bigr)
$$
where $\tau=\tau_0/|\theta|$.
The test rejects $H_0$ when $(\hat\theta/\tau_0)^2\ge \chi^{2,1-\alpha}_{(1)}$
and so it has power
\begin{align}\label{eq:thepower}
\Pr( \chi^{\prime,2}(\tau^{-2})\ge \chi^{2,1-\alpha}_{(1)}.
\end{align}
Given $\alpha$ and a level for power, we can solve~\eqref{eq:thepower}
for a value of $\tau$.
Some R code to do this is in the appendix.
Figure~\ref{fig:pow06} shows the results for
$\alpha=0.05$ and power $0.06$.

\begin{figure}
\centering
\includegraphics[width=.95\hsize]{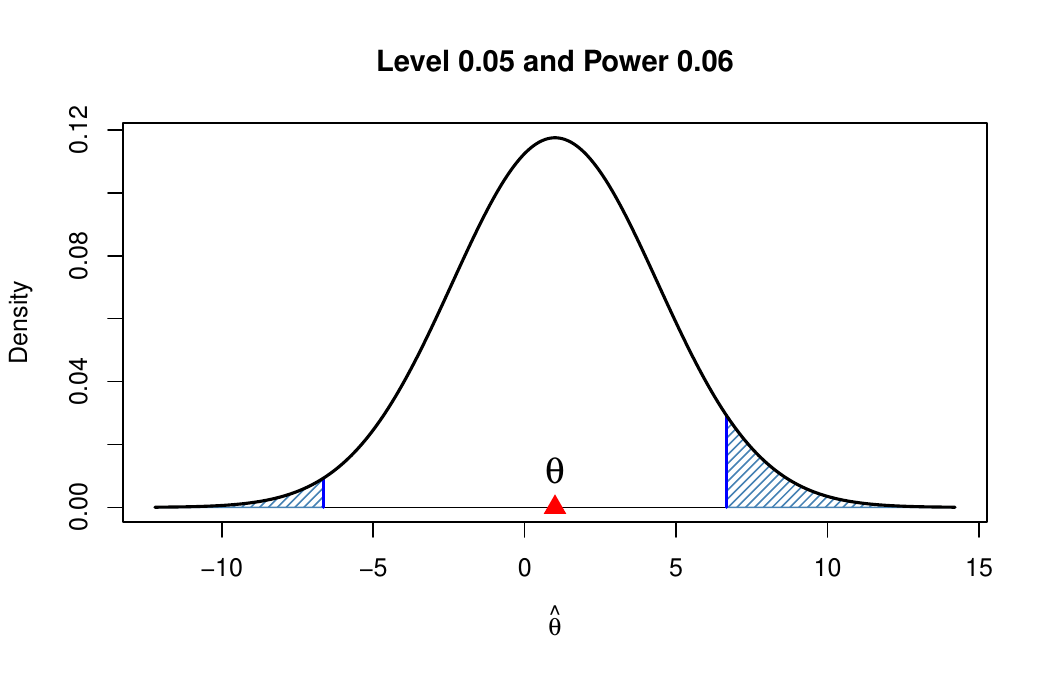}
\caption{\label{fig:pow06}
The curve shows the $\dnorm(1,\tau^2)$ density
of $\hat\theta$ when a test of $H_0:\theta=0$
at level $0.05$  has power $0.06$.
A triangle marks the true $\theta$.
The shaded regions where $|\hat\theta|\ge 6.65$
comprise the rejection region for $H_0$.
About $20$\% of rejections have the wrong sign.
}
\end{figure}

For $\alpha=0.05$ and power $0.06$, solving~\eqref{eq:thepower}
leads to $\tau\doteq3.39$.
Then the test rejects $H_0$ when
$|Z|\ge\Phi^{-1}(0.975)\times\tau \doteq 6.65$.  Therefore any significant
discovery must overestimate $|\theta|$ by at least $6.65$-fold.
When $\alpha=0.05$ and the power is $0.06$, the
average  value of $|\hat\theta/\theta|$ given that $H_0$
is rejected is $8.01$ and the
probability of a sign error is about $20$\%.
These quantities were computed by numerical quadrature
using R functions in the Appendix.

The exaggeration values and sign error
estimates in~\cite{gelm:2014:power06} (as of August 2019) 
do not match these results.  There, the wrong sign probability is
given as 24\% and the  minimum exaggeration factor of 9
is much higher than the computed value of $6.65$.
It appears that the figure in the blog is based on slightly
different numbers, corresponding to power less than 6\%.
For $\alpha=0.05$ and $1-\beta = 0.06$ we
find that $\tau=3.394507$. The point of all those digits in $\tau$
is that we can plug this value in to the {\tt retrodesign} R function from
\cite{gelm:carl:2014} and verify that it yields power $0.06$ when $\alpha=0.05$.
See Figure~\ref{fig:retrodesign}.  The numerical values
from~\cite{gelm:2014:power06}
appear to describe power between $5.54$\% and $5.55$\% because within that
range the minimum exaggeration factor is about 9.

\begin{figure}[t!]\centering
\begin{verbatim}
> retrodesign
function(A, s, alpha=.05, df=Inf, n.sims=10000){
z <- qt(1-alpha/2, df)
p.hi <- 1 - pt(z-A/s, df)
p.lo <- pt(-z-A/s, df)
power <- p.hi + p.lo
typeS <- p.lo/power
estimate <- A + s*rt(n.sims,df)
significant <- abs(estimate) > s*z
exaggeration <- mean(abs(estimate)[significant])/A
return(list(power=power, typeS=typeS, exaggeration=exaggeration))
}

> set.seed(1);retrodesign(1,3.394507)
$power
[1] 0.06

$typeS
[1] 0.2013426

$exaggeration
[1] 7.978919
\end{verbatim}
\caption{\label{fig:retrodesign}
R code from \cite{gelm:carl:2014} to verify that
standard error $\tau = 3.394507$ and effect size
$\theta=1$ yield power $0.06$ at $\alpha=0.05$.
The random seed is set to make the results more reproducible.
}
\end{figure}

\section{Confidence interval separation}\label{sec:ciwidth}

When evaluating a hypothesis test from the literature
we cannot tell what the power was.  There is often
an accompanying confidence interval.
Confidence intervals have the benefit of being in the same units as
the effect $\theta$ and hence they facilitate the study of practical significance.
Here they have the additional benefit that we can compare the
width of the confidence interval to its distance from the origin
to get an idea of how reliable $\sign(\hat\theta)$ is.

If the confidence interval is well separated from the
origin, then its creators could have widened it to a higher level
of confidence and still have excluded the origin.
Equivalently, they could have rejected $H_0$ at an
even smaller $\alpha$ than the one they used.

To illustrate, suppose that $H_0$ has been
rejected at level $\alpha=0.05$ because
$|\hat\theta|\ge \Phi^{-1}(0.975) s\doteq1.96s$.
Call this event $R$.
Suppose next,  that we only declare the sign of $\theta$
to match that of $\hat\theta$ when
the center $\hat\theta$ of the confidence interval is at
least $2\times1.96s$ away from $0$.  Call this event $S$.
That is, $S$ happens when the separation between the confidence
interval and $0$ is at least half of the confidence interval's width.
A sign error corresponds to $\hat\theta < -2\times1.96s$.
We call this event~$T$.

We will find $\Pr(S,T\giv R)$ in case $\theta>0$.
The probability is the same for $\theta<0$.
The conditional probability that we wrongly call the sign is
\begin{align*}
\Pr( S,T\giv R)
= \Pr( S,T,R)/\Pr(R) = \Pr(S,T)/\Pr(R)
\end{align*}
because $R$ occurs whenever $S$ does. Next
\begin{align*}
\Pr( S,T\giv R) & = \frac{\Pr( \hat\theta-\theta< -2\times 1.96s-\theta)}{\Pr( R)}\\
& \le \frac{\Pr( \hat\theta-\theta< -2\times 1.96s)}{\Pr( R)}\\
& = \frac{\Phi(-3.92)}{2\Phi(-1.96)}\\
& \doteq 0.0009.
\end{align*}
As before, we have used $s=\tau$, which is reasonable
when $\Pr( |s/\tau-1|>\epsilon)\to0$.

More generally,
suppose that $R$ is the event that a $100(1-\alpha)$\%
confidence interval for $\theta$ excludes the origin,
let $S$ be the event that $|\hat\theta-\theta|\ge\lambda\Phi^{-1}(1-\alpha/2)$
for some $\lambda\ge1$ and, assuming that $\theta\ne0$,
let $T$ be the event that $\sign(\hat\theta)\ne\sign(\theta)$.
Now suppose that we want
$\Pr( S,T\giv R)\le\alpha_S$
for a sign error quantity $\alpha_S$.
As in the example above,
\begin{align}\label{eq:pstrbound}
\Pr(S,T\giv R) \le  \frac{\Phi(-\lambda\Phi^{-1}(1-\alpha/2))}
{\alpha}.
\end{align}
We can keep the right side of~\eqref{eq:pstrbound}
below $\alpha_S$  by taking
$$
\lambda \ge
\frac{\Phi^{-1}(\alpha\alpha_S)}{\Phi^{-1}(\alpha/2)}=
\frac{\Phi^{-1}(1-\alpha\alpha_S)}{\Phi^{-1}(1-\alpha/2)}.
$$
The second expression above is a ratio of positive quantities.
Given an interval selected because $|\hat\theta|>\Phi^{-1}(1-\alpha/2)s$,
if we declare that $\sign(\theta)=\sign(\hat\theta)$ when
$|\hat\theta|>\lambda\Phi^{-1}(1-\alpha/2)s$, then the fraction
of such intervals where we declare the wrong sign will be at most $\alpha_S$.

We can describe this two step test in another way.  The first
step rejects $H_0$ if and only if $|\hat\theta|>\Phi^{-1}(1-\alpha_1/2)s$
for $\alpha_1\in(0,1)$.
The second step rejects $H_0$ if and only if $|\hat\theta|>\Phi^{-1}(1-\alpha_2/2)s$
for $\alpha_2\in(0,\alpha_1]$.
If both reject, then we have conditional probability at most
$$\frac{\alpha_2/2}{\alpha_1}=\frac{\alpha_2}{2\alpha_1}
$$
of making a sign error. The $2$ in the denominator
arises because at most half of the $\alpha_2$ rejections have the wrong sign.
The refiltering  described here is an extremely basic form of inference
after model selection, as described by \cite{fith:sun:tayl:2014}
who also give a comprehensive bibliography.
After confidence intervals are selected based on excluding the
origin, we use a second test to draw an inference on $\sign(\theta)$.

\begin{figure}
\includegraphics[width=0.95\hsize]{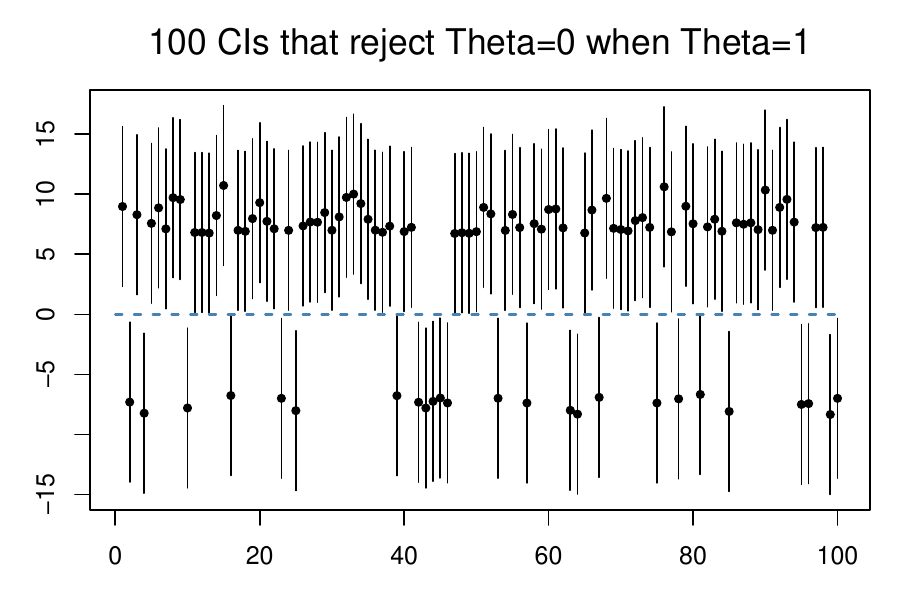}
\caption{\label{fig:cirej}
The figure shows 100 randomly selected confidence intervals
out of $10^5$ that were generated to have 6\% power at level $\alpha=0.05$.
The ones shown are the first 100 to exclude $0$.
}
\end{figure}

This approach leaves an awkward intermediate possibility
where we have rejected $H_0:\theta=0$ without
deciding whether $\theta$ is positive or negative.
This seems odd, but as we see next, it is reasonable. If we are only presented with
a confidence interval for $\theta$ that just barely
misses the origin, then we cannot be confident of $\sign(\theta)$.
Figure~\ref{fig:cirej} displays $100$ confidence intervals
selected to omit the origin, in a setting where the
test had 6\% power at the 5\% level.  These intervals
barely exclude the origin and many of them have center points
with the wrong sign.

\section{Some operating characteristics}\label{sec:power}

Here we make some numerical illustrations of how
the tests described here behave. We illustrate
with $\alpha = 0.05$, because that is a well known
level. Whether it is a good choice in practice depends
on context.  This  has been the subject of much
discussion \citep{wasserstein2016asa}.

\begin{figure}
\includegraphics[width=0.95\hsize]{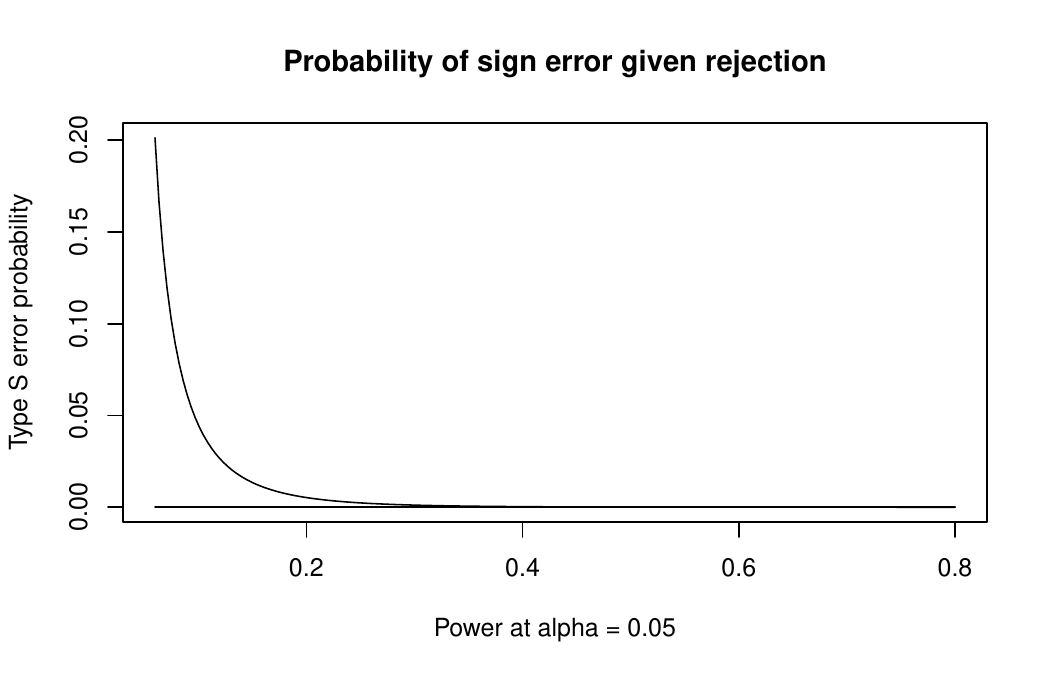}
\caption{\label{fig:types}
The figure shows how the probability of a sign error,
conditional on rejecting $H_0$, depends on power.
}
\end{figure}

Figure~\ref{fig:types} shows the conditional probability
of a sign error given that $H_0$ is rejected at the
5\% level.  The power of the underlying test ranges from 6\%
to 80\%.
The probability of a sign error becomes negligible at power greater than about 30\%.
\cite{gelm:carl:2014} also plot this quantity.

\begin{figure}
\includegraphics[width=0.95\hsize]{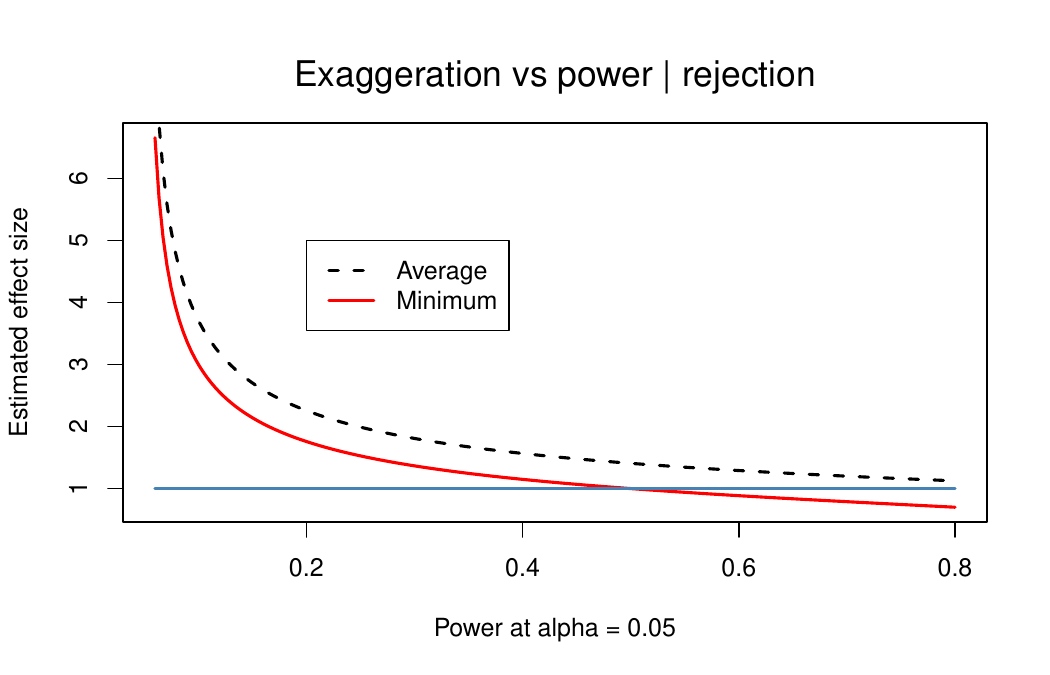}
\caption{\label{fig:exag}
The horizontal axis is the power of a test at level $\alpha=0.05$.
The solid curve is the expected value of $|\hat\theta/\theta|$
given that $H_0$ has been rejected.  The lighter curve is the
minimum value of $|\hat\theta/\theta|$ for which
$H_0$ could be rejected.
}
\end{figure}

Next we consider magnitude errors in $\hat\theta$.
The solid curve in Figure~\ref{fig:exag} shows the minimal value
of $|\hat\theta/\theta|$ for which $H_0$ can be rejected
at level $0.05$.  Whenever the power is below $0.5$
this minimal ratio is above $1$. When the power exceeds 50\%
then $H_0$ can be rejected with either an exaggerated magnitude
or an underestimated magnitude.  The  dashed curve in
Figure~\ref{fig:exag} shows the expected value of $|\hat\theta/\theta|$
given that $H_0$ has been rejected, versus the power of the test.
Even at 50\% power the average exaggeration can be as high as $1.4$ fold.
\cite{gelm:carl:2014} also plot the average exaggeration.

\begin{figure}
\includegraphics[width=0.95\hsize]{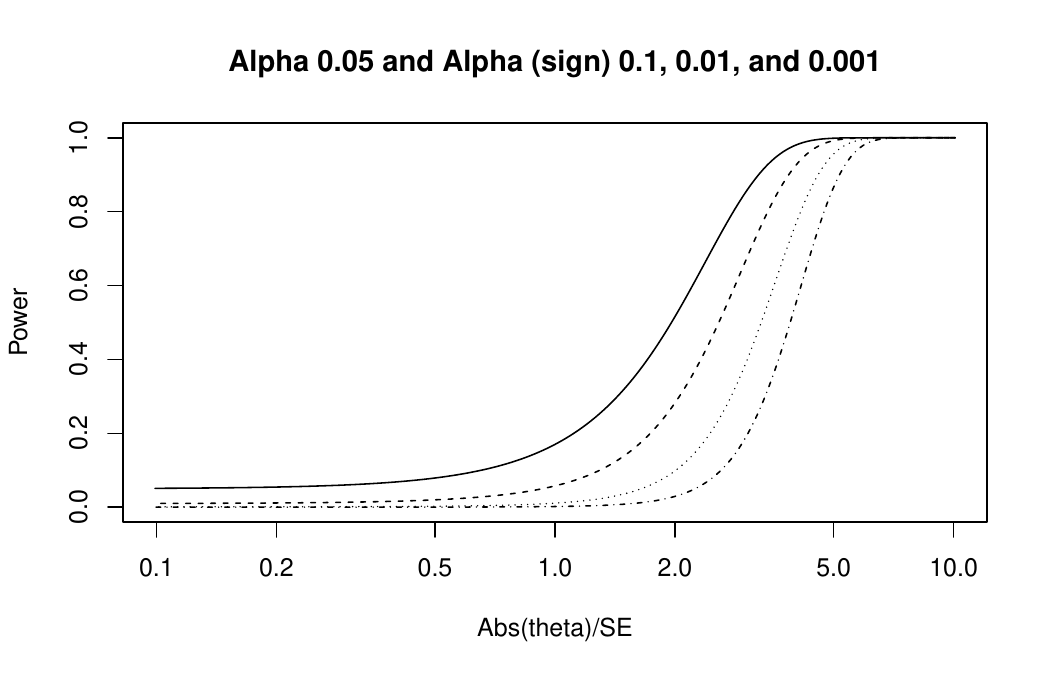}
\caption{\label{fig:powpow}
The horizontal axis is the effect size $|\theta|/\tau$.
The solid curve is the power of a test of
$H_0:\theta=0$ at level $\alpha =0.05$.
The next curves from top to bottom are the probabilities
of also making a definitive sign declaration
at $\alpha_S\in\{0.1,0.01,0.001\}$.
}
\end{figure}

Suppose for example, that we want to test $H_0$
at level $\alpha = 0.05$ but we are very averse
to sign errors and want $\alpha_S=0.001$.  That is,
at most one in one thousand rejections of $H_0$
should come with a sign wrongly declared by
the test at level $\alpha_2 = 0.0001$.
Figure~\ref{fig:powpow} shows both the power to reject $H_0$
and the power to identify $\sign(\theta)$ as a function of
$|\theta|/\tau$. It also depicts more lenient sign conditions
given by $\alpha_S=0.01$ and $\alpha_S=0.1$.
Ordinarily $\tau \doteq\sigma/\sqrt{n}$
where $\sigma$ is an asymptotic standard deviation.
There is a range of effect sizes in which rejections of $H_0$
without decisions on $\sign(\theta)$ will be common.
But once $|\theta|\sqrt{n}\ge 5\sigma$ or so, there is
very high power that $H_0$ will be rejected and a sign
will be identified.  The hypotheses in limbo
with $H_0$ rejected but no sign determination are
the ones for which obtaining additional data may be most valuable.

\section{Sign $p$-values}\label {sec:signp}
In high throughput settings, one might make a very
large number of primary hypothesis tests at level $\alpha_1$.
Suppose that they have been selected in a naive way, simply
providing the confidence intervals that omit zero and
withholding the others.
For each one that is rejected, there is then a smallest $\alpha_2\le\alpha_1$
at which $H_0$ is also rejected.  Call this value $p_2$.
Then the smallest $\alpha_S$ for which we would have
found the sign significant is $p_S = p_2/(2\alpha_1)$.  Note
that $p_S$ depends on the chosen value $\alpha_1$, which
we assume has been fixed in advance and not modified
to reduce or increase $p_S$.  Though we call it a $p$-value
because it governs a probability of misinterpretation, note
that it cannot be larger than~$1/2$.

\cite{benj:yeku:2005}
consider \cite{giov:asch:etal:1995}
who tested relationships between diet variables and cancers.
Their abstract reports three 95\% confidence intervals for relative risk.
The most promising one has a relative
risk estimate of $0.65$ with a 95\% confidence
interval from $0.44$ to $0.95$.  The face value
interpretation is of a protective effect versus
prostate cancer from consumption of tomatoes,
tomato sauce, tomato juice and pizza.
It is natural to measure relative risk on a logarithmic
scale. Their confidence interval is then
$-0.431\pm 0.379$, which is approximate due to rounding.
The 95\% confidence interval for lycopene intake
(which overlaps tomatoes et cetera) was,
on a log scale,
$-0.236\pm0.226$.
The third interval was approximately
$-0.755\pm0.755$.

We use these numerical examples below. We do not
make any claim about the relationship between diet and cancer.
That would require revisiting their data analyses, model
selections and considering the related literature on diet
and cancer.  Here we just show the values of $p_S$
corresponding to those intervals, in Table~\ref{tab:directionalp}.
None of those selected intervals would lead us to be
confident about $\sign(\theta)$.

\begin{table}
\centering
\begin{tabular}{ccccc}
\toprule
$\hat\theta$ & $1.96s$ &$\lambda=|\hat\theta|/1.96s$ &
$p_S=\Phi(-1.96\lambda)/0.05$\\
\midrule
$-.431$ & $0.379$ & $1.137$ & $0.258$\\
$-.236$ & $0.226$ & $1.044$ & $0.407$\\
$-.755$ & $0.755$ & $1.000$ & $0.500$\\
\bottomrule
\end{tabular}
\caption{\label{tab:directionalp}
This table shows some examples of $p_S$, the refiltered
$p$-value for $\sign(\theta)$.
}
\end{table}

\section{One-sided tests}\label{sec:onetail}

Here we investigate the possibility of sign errors in one-tailed tests.
Those tests are sometimes justified by an assumption
that the direction of the effect is certain a priori.
Then sign errors are not possible, unless our
certain opinion is wrong.  The a priori
uncertainty about the sign of $\theta$ must however be small
compared to the critical level $\alpha$ in use, and that is a very strong
assumption to work under.

A second justification is that sometimes only one direction
is consequential.  An opinion about one direction being
inconsequential might also be mistaken, but for sake of argument
we work with it. In this situation we might well make
a sign error by rejecting $H_0$.
Suppose that $\hat\theta\sim\dnorm(1,\tau^2)$,
so that $\hat\theta/\tau\sim\dnorm(1/\tau,1)$. The
true effect is positive but we might be testing for
a (consequential) negative effect.

A one-tailed test at level $\alpha$ in the negative direction would reject $H_0$
if $\hat\theta/\tau\le \Phi^{-1}(\alpha)$.  This happens with probability
$\Phi( \Phi^{-1}(\alpha)-1/\tau)$.

Now we revisit the case of $6$\% power for a two-tailed test at $\alpha = 0.05$.
Then $\tau \doteq 3.39$ and the sign error is $0.20$.
The wrong direction one-tailed test at level $0.05$ will reject $H_0$
with probability $\Phi( \Phi^{-1}(0.05)-1/\tau) \doteq 0.026$.

The chance of a wrong sign rejection is actually larger than $\alpha/2$ here.
That is not to say the conditional probability of a rejection being wrong
is over $50$\%.  Indeed in this setting any rejection at all is a sign error.
What it does mean is that in this low power setting, the probability of a sign error is
not small compared to $\alpha$. As a result, using one-tailed tests does
not correct the problem of sign errors in low power settings.

\section{Conclusions}\label{sec:conc}

The practice of reporting only the hypothesis
tests that are significant
at some level $\alpha$ is not recommended because
it can bring large sign and magnitude errors.
When one sees such results based on an asymptotically
normal test statistic, it is possible to filter them in such
a way that the tests passing the filter have a small probability
of a sign error.
This article presented one such method based on refiltering.
It can be applied without assuming
anything about the dependencies
among  the corresponding test statistics.

\section*{Acknowledgments}
I thank Jelena Markovic and Andrew Gelman for discussions
on an earlier version of this article.  Thanks to
Daniel Yekutieli, Aaditya Ramdas, and Edgar Dobriban
for further discussions at WHOA-PSI 4.
This work was supported by the U.S.\ National Science
Foundation under
grants DMS-1407397,  DMS-1521145 and IIS-1837931.

\bibliographystyle{apalike}
\bibliography{powersix}

\vfill\eject
\section*{Appendix: some R code}

This function solves equation~\eqref{eq:thepower} for $\tau$
given the level and power of a test.
\begin{verbatim}
tau = function(alpha,powr){
# A Gaussian test at level alpha
# has power powr for Y ~ N( 1, tau^2 ).
# Solve for tau using noncentral chisquare.

aux = function(tau){
  1-powr - pchisq(qchisq(1-alpha,1),1,ncp=1/tau^2)
}

ur = uniroot( aux, lower=10^-9,upper=10^6,
    tol = .Machine$double.eps^0.9)
ur$root
}
\end{verbatim}

The function {\tt exag} below computes the
expected value of $\hat\theta$ conditionally
on $H_0$ being rejected and on $\sign(\hat\theta)$.
It uses {\tt gauscondmean} which
computes $\e(Z\giv A\le Z\le B)$ for $Z\sim\dnorm(\mu,\sigma^2)$.

\begin{verbatim}
exag = function(alpha=.05,powr=.06){
# Compute exaggeration factors

critse  = tau(alpha,powr)  # Critical Standard Error
critz   = critse*qnorm(1-alpha/2)
posmean = gauscondmean(mu=1,sigma=critse,A=critz, B=Inf, n=10^5 )
negmean = gauscondmean(mu=1,sigma=critse,B=-critz, n=10^5 )
typeS   = pnorm(-critz,1,critse)/powr

exager = posmean*(1-typeS)+abs(negmean)*typeS

list(critse=critse,critz=critz,exager=exager,posmean=posmean,negmean=negmean,typeS=typeS)
}
\end{verbatim}

The function {\tt gauscondmean} below
computes $\e(Z\giv A\le Z\le B)$ for $Z\sim\dnorm(\mu,\sigma^2)$
using {\tt stdgauscondmean} which handles the case $\mu=0$ and $\sigma=1$.

\begin{verbatim}
gauscondmean = function(mu=0,sigma=1,A=-Inf, B=Inf, n=10^5 ){
# Average of N( mu, sigma^2 ) over interval (A,B)

mu + abs(sigma)*stdgauscondmean((A-mu)/abs(sigma),(B-mu)/abs(sigma),n)
}
\end{verbatim}

The function {\tt stdgauscondmean} below
computes $\e(Z\giv A\le Z\le B)$ for $Z\sim\dnorm(0,1)$.
It uses a midpoint rule.

\begin{verbatim}
stdgauscondmean = function(A=-Inf, B=Inf, n=10^5, plot=FALSE ){
# Average of N( 0, 1 ) over interval (A,B)

u = ((1:n)-0.5)/n
u = pnorm(A) + (pnorm(B)-pnorm(A)) * u
z = qnorm(u)
if( plot )hist(z,50) # for testing/debug
mean(z)
}
\end{verbatim}

\end{document}